\documentclass[a4paper,12pt]{article}
\usepackage{a4wide}
\usepackage{amsmath}
\usepackage{amssymb}
\usepackage{amsthm}
\usepackage{latexsym}
\usepackage{graphicx}
\usepackage[english]{babel}
\usepackage{makeidx}

\newtheorem{conj}[subsection]{Conjecture}
\newtheorem{teor}[subsection]{Theorem}

\newtheorem{cor} [subsection]{Corollary}
\newcommand{\Zng}{$\mathbb Z^n$-graded $S$-module}
\newcommand{\me}{\mathbf m}
\def\sdepth{\operatorname{sdepth}}
\def\depth{\operatorname{depth}}
\begin{document}
\selectlanguage{english}
\frenchspacing

\large
\begin{center}
\textbf{Stanley depth of square free Veronese ideals}

by

Mircea Cimpoea\c s
\end{center}
\normalsize

\begin{abstract}
We compute the Stanley depth for the quotient ring of a square free Veronese ideal and we give some bounds for
the Stanley depth of a square free Veronese ideal. In particular, it follows that both satisfy the Stanley's conjecture.


\end{abstract}

\section*{Introduction}

Let $K$ be a field and $S=K[x_1,\ldots,x_n]$ the polynomial ring over $K$.
Let $M$ be a \Zng. A \emph{Stanley decomposition} of $M$ is a direct sum $\mathcal D: M = \bigoplus_{i=1}^rm_i K[Z_i]$ as $K$-vector space, where $m_i\in M$, $Z_i\subset\{x_1,\ldots,x_n\}$ such that $m_i K[Z_i]$ is a free $K[Z_i]$-module. We define $\sdepth(\mathcal D)=min_{i=1}^r |Z_i|$ and $\sdepth(M)=max\{\sdepth(M)|\;\mathcal D$ is a Stanley decomposition of $M\}$. The number $\sdepth(M)$ is called the \emph{Stanley depth} of $M$. Herzog, Vladoiu and Zheng show in \cite{hvz} that this invariant can be computed in a finite number of steps if $M=I/J$, where $J\subset I\subset S$ are monomial ideals.

There are two important particular cases. If $I\subset S$ is a monomial ideal, we are interested in computing $\sdepth(S/I)$ and $\sdepth(I)$. There are some papers regarding this problem, like \cite{hvz},\cite{asia},\cite{sum}, \cite{shen} and \cite{mir}. Stanley's conjecture says that $\sdepth(S/I)\geq \depth(S/I)$, or in the general case, $\sdepth(M)\geq \depth(M)$, where $M$ is a finitely generated multigraded $S$-module. The Stanley conjecture for $S/I$ was proved for $n\leq 5$ and in other special cases, but it remains open in the general case. See for instance, \cite{apel}, \cite{hsy}, \cite{jah}, \cite{pops}, \cite{popi} and \cite{pope}.



For any $d\in [n]$, we denote $I_{n,d}:=(u\in S$ square free monomial $:\; deg(u)=d)$. It is well known that $dim(S/I_{n,d})=\depth(S/I_{n,d})=d-1$. We show that $\sdepth(S/I_{n,d})=d-1$ and we give some bounds for $\sdepth(I_{n,d})$, see Theorem $1.1$. As a consequence, it follows that $I_{n,d}$ and $S/I_{n,d}$ satisfy the Stanley's conjecture, see Corollary $1.2$. Also, we prove that $\sdepth(I_{n,d})=d+1$, if $2d+1 \leq n \leq 3d$, see Corollary $1.5$. In order to do so, we prove some combinatorics results, see Theorem $1.3$ and Corollary $1.4$.

We conjecture that $\sdepth(I_{n,d})=\frac{n-d}{d+1} + d$, see $1.6$.

\footnotetext[1]{This paper was supported by CNCSIS, ID-PCE, 51/2007}

\newpage
\section{Main results}

\begin{teor}
(1) $\sdepth(S/I_{n,d})=d-1$.

(2) $d\leq \sdepth(I_{n,d}) \leq \frac{n-d}{d+1} + d$.
\end{teor}

\begin{proof}
(1) Firstly, note that $\sdepth(S/I_{n,d})\leq d-1 = dim(S/I_{n,d})$. We use induction on $n$ and $d$. If $n=1$, there is nothing to prove. If $d=1$, it follows that $I_{n,1}=(x_1,\ldots,x_n)$ and thus $\sdepth(S/I_{n,1})=0$, as required. If $d=n$, it follows that $I_{n,n}=(x_1\cdots x_n)$ and therefore $\sdepth(S/I_{n,n})=n-1$, as required. Now, assume $n>1$ and $1<d<n$. Note that 
\[ S/I_{n,d}=\bigoplus_{|supp(u)|<d} u\cdot K = \sum_{Z\subset \{x_1,\ldots,x_n\},\;|Z|=d-1}K[Z].\]
We denote $S'=K[x_1,\ldots,x_{n-1}]$. By previous equality, we get $ S/I_{n,d} =$ \linebreak \[ \sum_{Z\subset \{x_1,\ldots,x_{n-1}\},\;|Z|=d-1}K[Z] 
\oplus x_n (\sum_{Z\subset\{x_1,\ldots,x_{n-1}\},\;|Z|=d-1}K[Z])[x_n] = S'/I_{n-1,d} \oplus x_n(S'/I_{n-1,d-1})[x_n].\]
By induction hypothesis, it follows that $\sdepth(S/I_{n,d})=d-1$. 

(2) We consider the following simplicial complex, associated to $I_{n,d}$, 
\[ \Delta_{n,d}:= \{supp(u):\; u\in I_{n,d} \;monomial\}. \]
Note that, by \cite[Theorem 2.4]{hvz}, there exists a partition of $\Delta_{n,d}=\bigcup_{i=1}^r [F_i,G_i]$, such that $\min_{i=1}^r |G_i| = \sdepth(I_{n,d}):=s$. Note that $\Delta_{n,d}=\{F\subset [n]:\; |F|\geq d \}$. It follows that $\sdepth(I_{n,d})\geq d$.

We consider an interval $[F_i,G_i]$ with $|F_i|=d$. Since $|G_i|\geq s$, it follows that there exists at least $(s-d)$ distinct sets in $[F_i,G_i]$ of cardinality $d+1$. Since $\Delta_{n,d} = \bigcup_{i=1}^r [F_i,G_i]$ is a partition, it follows that
\[ \binom{n}{d+1} = \frac{n-d}{d+1}\binom{n}{d}  \geq (s - d) \binom{n}{d}. \]
Thus, $s\leq d + \frac{n-d}{d+1}$.
\end{proof}

\begin{cor}
$I_{n,d}$ and $S/I_{n,d}$ satisfy the Stanley's conjecture. Also, \[ \sdepth(I_{n,d})\geq \sdepth(S/I_{n,d})+1.\]
\end{cor}

Let $k\leq n$ be two positive integers. We denote $A_{n,k}=\{F\subset [n]|\; |F|=k\}$.

\begin{teor}
For any positive integers $d\leq n$ such that $d\leq n/2$, there exists a bijective map $\Phi_{n,d}:A_{n,d} \rightarrow A_{n,d}$, such that $\Phi_{n,d}(F)\cap F = \emptyset$ for any $F\in A_{n,d}$.
\end{teor}

\begin{proof}
We use induction on $n$ and $d$. If $n\leq 2$ the statement is obvious. If $d=1$, for any $i\in [n]$, we define $\Phi_{n,1}(\{i\})=\{j\}$, where $j=\max ([n]\setminus \{\Phi_{n,1}(\{1\}),\ldots,\Phi_{n,1}(\{i-1\})\})$. $\Phi_{n,1}$ is well defined and satisfy the required conditions.
 
Now, assume $n\geq 3$ and $d\geq 2$. If $n=2d$ we define $\Phi_{n,d}(F)=[n]\setminus F$. Obviously, $\Phi_{n,d}$ satisfy the required conditions. Thus, we may also assume $d<n/2$.

On $A_{n,d}$, we consider the lexicographic order, recursively defined by $F<G$ if and only if $max\{F\}<max\{G\}$ or $max\{F\} = max\{G\} = k$ and $F\setminus \{k\} < G\setminus \{k\}$ on $A_{n,d-1}$. For any $F\in A_{n,d}$, we define $G:=\Phi_{n,d}(F)$ to be the maximum set, with respect to "$<$", such that $G \cap F=\emptyset$ and $G \neq \Phi_{n,d}(H)$ for all $H<F$. In order to complete the proof, it is enough to show that each collection of sets \[ \mathcal M^n_F = \{ G\subset[n]\;:\; |G|=d,\; G\cap F=\emptyset, \; G \neq \Phi_{n,d}(H)\;\; (\forall)\; H<F \}\] is nonempty, for all $F\subset [n]$. Assume there exists some $F\subset [n-1]$ such that $\mathcal M^n_F=\emptyset$. It obviously follows that $M^{n-1}_F=\emptyset$ and thus $\Phi_{n-1,d}$ is not well defined, a contradiction. Also, if $M_F=\emptyset$ for some $F\subset [n]$ with $n\in F$, it follows similarly that $\Phi_{n-1,d-1}$ is not well defined, again a contradiction. Therefore, the required conclusion follows.
\end{proof}

\begin{cor}
For any positive integers $d$ and $n$ such that $d < n/2$, there exists an injective map $\Psi_{n,d}:A_{n,d} \rightarrow A_{n,d+1}$, such that $F\subset \Psi_{n,d}(F)$ for any $F\in A_{n,d}$.
\end{cor}

\begin{proof}
We use induction on $n$. If $n\leq 2$ there is nothing to prove. If $d=1$, we define $\Psi_{n,1}:A_{n,1}\rightarrow A_{n,2}$ by $\Psi_{n,1}(\{1\})=\{1,2\}$, \ldots, $\Psi_{n,1}(\{n-1\})=\{n-1,n\}$ and $\Psi_{n,1}(\{n\})=\{1,n\}$. Now, assume $n\geq 3$ and $d\geq 2$.

If $n=2d+1$, we consider the bijective map $\Phi_{n,d}: A_{n,d} \rightarrow A_{n,d}$ such that $\phi(F)\cap F=\emptyset$ for all $F\in A_{n,d}$ and we define $\Psi_{n,d}(F):=[n]\setminus \Phi_{n,d}(F)$. The map $\Psi_{n,d}$ satisfies the required condition.

If $n<2d+1$, we define $\Psi_{n,d}(F):= \Psi_{n-1,d}(F)$ if $F\subset [n-1]$ and $\Psi_{n,d}(F):= \Psi_{n-1,d-1}(F\setminus \{n\}) \cup \{n\}$ if $n\in F$. Note that both $\Psi_{n-1,d}$ and $\Psi_{n-1,d-1}$ are well defined and injective by induction hypothesis, since $n-1\leq 2d+1$. It follows that $\Psi_{n,d}$ is well defined and injective, as required.
\end{proof}

\begin{cor}
Let $n,d$ be two positive integers such that $2d+1 \leq n \leq 3d$. Then $\sdepth(I_{n,d})=d+1$.
\end{cor}

\begin{proof}
As in the proof of $1.1$, we denote $\Delta_{n,d}:= \{supp(u):\; u\in I_{n,d} \;monomial\}=\{F\subset [n]:\; |F|\geq d \}$. We consider the following partition of $\Delta_{n,d}$:
\[ \Delta_{n,d} = \bigcup_{|F|=d}[F,\Psi_{n,d}(F)] \cup \bigcup_{|F|>d+1}[F,F], \]
where $\Psi_{n,d}$ is given by the previous corollary. It follows that $\sdepth(I_{n,d})\geq d+1$. On the other hand, by $1.1$, $\sdepth(I_{n,d})\leq d+1$ and thus $\sdepth(I_{n,d}) = d+1$, as required.
\end{proof}

\begin{conj}
For any positive integers $d\leq n$ such that $d\leq n/2$, $\sdepth(I_{n,d})=\frac{n-d}{d+1} + d$
\end{conj}

\vspace{2mm} \noindent {\footnotesize
\begin{minipage}[b]{15cm}
 Mircea Cimpoeas, Institute of Mathematics of the Romanian Academy, Bucharest, Romania\\
 E-mail: mircea.cimpoeas@imar.ro
\end{minipage}}
\end{document}